\documentclass{amsart}
\title%[Further extended Bailey's transforms and applications]
{Positive alternating sums of integer partitions}
\usepackage{amssymb,amsmath,amsthm,epsfig,graphics,latexsym,hyperref}
\theoremstyle{definition}

\theoremstyle{plain}
\newtheorem{lemma}      {Lemma}

\newtheorem{theorem}    {Theorem}

\newtheorem{corollary}  {Corollary}

\theoremstyle{remark}

\newcommand{\fr}{\frac}

 \newmuskip\pFqskip
\mathchardef\pFcomma=\mathcode`, % keep a copy of the comma

\begin{document}
  \author[M. El Bachraoui]{Mohamed El Bachraoui}
  \address{Dept. Math. Sci,
 United Arab Emirates University, PO Box 15551, Al-Ain, UAE}
 \email{melbachraoui@uaeu.ac.ae}
% \address{American University of Sharjah, Department of Mathematics and Statistics,
%American University of Sharjah, P.O. Box 26666, Sharjah, UAE}
%\email{jgriffin@aus.edu}
%
\begin{abstract}
We combine an extended version of Bailey's transform with an identity of Bressoud and with some
identities of Berkovich and Warnaar to prove a variety
of positivity results for alternating sums involving partition functions.
\end{abstract}
\date{\textit{\today}}
\maketitle
\section{Introduction}
Throughout the paper $n$ is a nonnegative integer.
The $q$-shifted factorials are given by
\[
(a;q)_0 = (a)_0 = 1,\ \quad, (a;q)_n = (a)_n = \prod_{j=0}^{n-1} (1-aq^j),\quad
(a;q)_{\infty} = (a)_{\infty} = \prod_{j=0}^{\infty} (1-aq^j),
\]
and the $q$-binomial coefficient is given for any nonnegative integers $M$ and $N$ by
\[
{N\brack M}_q = {N\brack M} =
\begin{cases}
\fr{(q)_N}{(q)_M (q)_{N-M}} &\text{if\ } N\geq M, \\
0 & \text{otherwise.}
\end{cases}
\]
It is well-known that for integers $N\geq M \geq 0$, the $q$-binomial coefficient
${N\brack M}$ is a polynomial of degree $N(N-M)$ which has nonnegative integer coefficients.
As usual we let $p(n)$ denote the number of partitions of $n$. Information about $q$-binomial coefficients
and integer partition can be found in the book by Andrews~\cite{Andrews}.
Following~\cite{Corteel-Lovejoy} an overpartition
of $n$ is a partition with the extra condition that the first occurrence of any part may be overlined or not
and the number of such overpartitions is written $\overline{p}(n)$. As to generating functions, we have
\begin{equation}\label{generating}
\sum_{n=0}^{\infty} p(n)q^n = \fr{1}{(q;q)_{\infty}} \ \text{and\ }
\sum_{n=0}^{\infty} \overline{p}(n)q^n = \fr{(-q;q)_{\infty}}{(q;q)_{\infty}}, \quad
(p(0)=\overline{p}(0)= 1).
\end{equation}
Alternating sums of integer partitions have received much attention by mathematicians in recent years.
For paper on this subject which are related to our work, see for instance
\cite{Andrews-Merca-2012, Andrews-Merca-2018, Guo-Zeng, Merca}.
Merca~\cite{Merca} proved that for all $n>0$
\[
p(n-5)-p(n-2)-p(n-1)+p(n) \geq 0
\]
which Andrews and Merca~\cite{Andrews-Merca-2012} generalized by showing that for $n, k >0$,
\begin{equation}\label{AndMer-ineq}
(-1)^{k-1}\sum_{j=0}^{k-1} (-1)^j\Big(p(n-j(3j+1)/2)-p(n-j(3j+5)/2-1) \Big) \geq 0.
\end{equation}
Their proof is partition-theoretic but it is based on the following truncated sum
\begin{equation}\label{AndMer-truncated-0}
\fr{1}{(q;q)_{\infty}} \sum_{j=0}^{k-1} (-1)^j q^{j(3j+1)/2}(1-q^{2j+1})
=
1+ (-1)^{k-1} \sum_{n=1}^{\infty}\fr{q^{{k\choose 2}+(k+1)n}}{(q;q)_n} {n-1 \brack k-1}.
\end{equation}
To derive the foregoing identity, Andrews and Merca used induction along with the $q$-binomial
theorem and the basic properties of the $q$-binomial coefficients.
Guo and Zeng~\cite{Guo-Zeng} established the following truncated sum
\begin{equation}\label{GuZe-truncated-1}
\fr{(-q)_{\infty}}{(q)_{\infty}} \Big(1+ 2\sum_{j=1}^{k} (-1)^j q^{j^2} \Big)
=
1+ (-1)^{k} \sum_{n=k+1}^{\infty}\fr{(-q)_k (-1)_{n-k} q^{(k+1)n}}{(q)_n} {n-1 \brack k-1}
\end{equation}
and as a consequence they obtained
\begin{equation}\label{GuZe-ineq-1}
(-1)^{k}\Big( \overline{p}(n) + 2\sum_{j=0}^{k} (-1)^j\overline{p}(n-j^2) \Big) \geq 0.
\end{equation}
The authors' proof for (\ref{GuZe-truncated-1}) relies on the following formula of
Andrews~\cite{Andrews 1986}
\[
\sum_{j=0}^n \fr{(b)_j (1-bq^{2j})(b/a)_j a^j q^{j^2}}{(1-b) (q)_j (aq)_j}
=
\fr{(bq)_n}{(aq)_n} \sum_{j=0}^n \fr{(b/a)_j a^j q^{(n+1)j}}{(q)_j}
\]
which is a limiting case of Watson $q$-Whipple transformation, see Gasper and Rahman~\cite[III. 18]{Gasper-Rahman} and Andrews~\cite[Lemma 2]{Andrews 1986}.
Besides, Guo and Zeng~\cite[p. 702]{Guo-Zeng} stated that it is an open problem to give
partition-theoretic interpretation for (\ref{GuZe-ineq-1}).
Recently, Andrews and Merca~\cite{Andrews-Merca-2018}
provided such interpretations by first rewriting (\ref{GuZe-truncated-1})
in the following form
\begin{equation}\label{AndMer-truncated-1}
\fr{(-q)_{\infty}}{(q)_{\infty}}
\Big( 1+2\sum_{j=1}^k (-1)^j q^{j^2} \Big)
\end{equation}
\[
=
1+ 2 (-1)^k \fr{(-q)_k}{(q)_k}
\sum_{j=0}^{\infty} \fr{q^{(k+1)(k+j+1)} (-q^{k+j+2})_{\infty}}
{(1-q^{j+k+1})(q^{j+k+2})_{\infty}}
\]
and then explaining the right-hand side of this formula in terms integer partitions.
The main argument in the authors' proofs for (\ref{AndMer-truncated-1})
is the following formula which is due to Rogers and Fine
(see \cite[p. 15]{Rogers})
\begin{equation}\label{Rogers-Fine}
\sum_{j=0}^{\infty}\fr{(a)_j c^j}{(b)_j} =\sum_{j=0}^{\infty}\fr{ (a)_j (acq/b)_j (bc)^j q^{j^2-j} (1-ac q^{2j})}{(b)_j (c)_{j+1}}.
\end{equation}
Furthermore, Andrews and Merca~\cite{Andrews-Merca-2018} showed that both
(\ref{AndMer-truncated-0}) and (\ref{GuZe-truncated-1}) can be established
as instances of  (\ref{Rogers-Fine}).
In this paper we shall establish positivity results for alternating sums of the partition functions
$p(n)$ and $\overline{p}(n)$. We do that by showing that their corresponding generating functions
have nonnegative coefficients.
Specifically, we will prove the following main results.
\begin{theorem}\label{thm result-1}
There holds
\[
\fr{1}{(q;q)_{\infty}} \sum_{n=0}^{\infty} (-1)^n q^{\fr{13n^2+11n}{2}}(1-q^{2n+1})
 =
\sum_{n=0}^{\infty} \fr{q^{n^2+2n}}{(q^2;q)_{2n+1}} \sum_{k=0}^n q^{k^2+k} {n\brack k}_q.
\]
\end{theorem}
\begin{corollary}\label{cor result-1}
\[
\sum_{j=0}^n (-1)^j \Big( p \big(n-j(13j+11)/2 \big)-p \big(n-j(13j+15)/2-1 \big) \Big) \geq 0.
\]
\end{corollary}
\begin{theorem}\label{thm result-2}
We have
\[
\fr{1}{(q;q)_{\infty}}\sum_{n=0}^{\infty} (-1)^n q^{4n (10n-1)+4} (1-q^{16n+8})
\]
\[
=
\sum_{n=0}^{\infty}\fr{(-q^{-1};q^2)_n}{(q;q)_{2n+3}}
\sum_{k=0}^n \fr{q^{(n-2k)^2+5n+2k}(-q^4;q^4)_{k+1}}{(1-q^{4k+2})(-q^{3-2n};q^2)_{2k}}
{2n \brack 2k}_{q^2}.
\]
\end{theorem}
\begin{corollary}\label{cor result-2}
\[
\sum_{j=0}^n (-1)^j \Big( \overline{p} \big(n-4j(10j+1)-4 \big)-
\overline{p} \big(n-4j(10j+3)-12 \big) \Big) \geq 0.
\]
\end{corollary}
\begin{theorem}\label{thm result-3}
We have
\[
\fr{1}{(q;q)_{\infty}} \sum_{n=0}^{\infty} (-1)^n q^{8n(3n+1)}(1-q^{32n+16})
\]
\[
= \sum_{n=0}^{\infty}\fr{(-q^{-1};q^2)_n}{(q^2;q^2)_n (q;q^2)_{n+2}}
\sum_{k= 0}^n \fr{q^{(n-2k)^2+5n+6k}(1+q^{4k+4})}{1-q^{4k+2}} {2n \brack 2k}_{q^2}.
\]
\end{theorem}
\begin{corollary}\label{cor result-3}
\[
\sum_{j=0}^n (-1)^j \Big( p \big(n-8j(3j+1) \big)-
p \big(n-8j(3j+5)-16 \big) \Big) \geq 0.
\]
\end{corollary}
\begin{theorem}\label{thm result-4}
There holds
\[
\fr{1}{(q;q)_{\infty}} \sum_{n=0}^{\infty} (-1)^n q^{6n(9n + 7)} (1-q^{24n+12})
\]
\[
=\sum_{n=0}^{\infty}\fr{(-q^{-1};q^2)_n}{(q^2;q^2)_n (q;q^2)_{n+3}}
\sum_{k=0}^n \fr{q^{(n-2k)^2+7n+2k}(1+q^4+q^8) (1+q^{6+4k}) }
{(1-q^{2k+2})(1-q^{2k+4}) (-q^{3-2n};q^2)_{2k}} {2n \brack 2k}_{q^2}.
\]
\end{theorem}
\begin{corollary}\label{cor result-4}
\[
\sum_{j=0}^n (-1)^j \Big( \overline{p} \big(n-6j(9j+7) \big)-
\overline{p} \big(n-6j(9j+11)-12 \big) \Big) \geq 0.
\]
\end{corollary}
\noindent
The key argument of our proofs is the following extension of the Bailey's transform.
\begin{lemma}\label{lem transf-1}\cite[Lemma 1]{Bachraoui}
Let $d, e, m$ and $n$ be positive integers such that $d\mid m$.
 Assuming convergence of the series, if
\[
B_n = \sum_{j=0}^{\lfloor{dn/m\rfloor}} A_{ej} U_{dn-mj} V_{dn+mj}
\quad\text{and\quad}
C_n = \sum_{j=mn/d}^{\infty} D_{ej} U_{dj-mn} V_{dj+mn},
\]
then
\[
\sum_{n=0}^{\infty} A_{en} C_n = \sum_{n=0}^{\infty} B_n D_{en}.
\]
\end{lemma}
\noindent
Note that the case $d=e=m=1$ in Lemma~\ref{lem transf-1}  gives the classical
Bailey's transform~\cite{Bailey-1}. For a survey on Bailey's transform along with applications we
refer for instance to
\cite{Andrews 2000, Andrews-Berkovich, Warnaar-2001, Warnaar-2003, Berkovich-Warnaar}.
We need the following well-known result
\begin{equation}\label{Bailey 0}
\sum_{k=n}^{\infty}\fr{a^k q^{k^2}}{(q)_{k-n} (aq)_{k+n}}
= \fr{a^{n} q^{n^2}}{(aq)_{\infty}},
\end{equation}
see for instance Bressoud~\cite[p. 215]{Bressoud}, Schilling and Warnaar~\cite[(1.3)]{Schilling-Warnaar},
and Warnaar~\cite[(1.5)]{Warnaar-2001b}.
Note that according to Schilling and Warnaar~\cite{Schilling-Warnaar}, the pair of sequences
$(\delta_k)_{k\geq 0}$ and $(\gamma_k)_{k\geq 0}$ where
\[
\delta_k = a^k q^{k^2} \quad \text{and\quad} \gamma_k = \fr{a^k q^{k^2}}{(aq)_{\infty}}
\]
is a \emph{conjugate Bailey pair} as they are connected through the relation~(\ref{Bailey 0}).
We refer the reader to~\cite{Schilling-Warnaar, Warnaar-2001b} for an survey on Bailey pairs, conjugate Bailey pairs, and applications.
\section{Proof of Theorem~\ref{thm result-1} and Corollary~\ref{cor result-1}}
\noindent
We need the following identity of Bressoud~\cite[(3.4)]{Bressoud}
\begin{equation}\label{Bressoud 3-4}
\sum_{k=0}^{n} \fr{a^k q^{k^2}{n\brack k}_q (a\gamma q;q)_n}
{(a\gamma q;q)_k (a\gamma q;q)_{n-k}}
\end{equation}
\[
= \fr{(a^2 q;q)_{2n}(q;q)_n}{(aq;q)_n} \sum_{k=0}^{\lfloor n/2 \rfloor}
\fr{a^{2k}\gamma^k q^{2k^2}(1-aq^{2k})(a;q)_k (\gamma^{-1};q)_k}
{(1-a) (q;q)_k (a\gamma q;q)_k (q;q)_{n-2k} (a^2 q;q)_{n+2k}}.
\]
We want to apply Lemma~\ref{lem transf-1} to (\ref{Bressoud 3-4}) as follows.
Let $m=2$, $d= e=1$ and let
\[
A_n= \fr{a^{2n}\gamma^n q^{2 n^2}(1-aq^{2n})(a;q)_n (\gamma^{-1};q)_n}{(1-a) (q;q)_n (a\gamma q;q)_n},\quad
D_n = a^{2n} q^{n^2},
\]
\[
U_n=\fr{1}{(q;q)_n},\quad V_n=\fr{1}{(a^2 q;q)_n}.
\]
Then
\[
\begin{split}
B_n &= \sum_{k=0}^{\lfloor n/2\rfloor} A_k U_{n-2k}V_{n+2k} \\
&= \sum_{k=0}^{\lfloor n/2\rfloor}\fr{a^{2k}\gamma^k q^{2k^2} (1-aq^{2k})(a;q)_k(\gamma^{-1};q)_k}
{(1-a)(q;q)_k (a\gamma q;q)_k (q;q)_{n-2k}(a^2 q;q)_{n+2k}} \\
&= \fr{(aq;q)_n}{(a^2 q;q)_{2n}(q;q)_n} \sum_{k=0}^n a^k q^{k^2} {n\brack k}_q \fr{(a\gamma q;q)_n}{(a\gamma q;q)_k (a\gamma q;q)_{n-k}}
\end{split}
\]
where the last identity follows by (\ref{Bressoud 3-4}).
Moreover, by virtue of identity~(\ref{Bailey 0}) we have
\[
C_n = \sum_{k=2n}^{\infty}D_k U_{k-2n}V_{k+2n}
= \sum_{k=2n}^{\infty} \fr{a^{2k}q^{k^2}}{(q;q)_{k-2n} (a^2 q;q)_{k+2n}}
= \fr{a^{4n} q^{4n^2}}{(a^2 q;q)_{\infty}}.
\]
Then by Lemma~\ref{lem transf-1}, we get
\begin{equation}\label{help-0}
\sum_{n=0}^{\infty} B_n D_n =
\sum_{n=0}^{\infty}\fr{a^{2n}q^{n^2}(aq;q)_n (a\gamma q;q)_n}{(a^2q;q)_{2n}(q;q)_n}
\sum_{k=0}^n\fr{a^k q^{k^2}{n\brack k}_q}{(a\gamma q;q)_k (a\gamma q;q)_{n-k}}
\end{equation}
\[
= \sum_{n=0}^{\infty} A_n C_n =
\fr{1}{(a^2 q;q)_{\infty}}
\sum_{n=0}^{\infty}\fr{a^{6n} \gamma^n q^{6n^2}(1-q^{2n})(a;q)_n (\gamma^{-1};q)_n}
{(1-a)  (q;q)_n (a\gamma q;q)_n}.
\]
Letting $\gamma \to 0$ in~(\ref{help-0}) yields
\begin{equation}\label{help-1}
\fr{1}{(a^2 q;q)_{\infty}} \sum_{n=0}^{\infty} \fr{(-1)^n a^{6n} q^{\fr{13n^2-n}{2}}(1-aq^{2n}) (a;q)_n}
{(1-a)(q;q)_n}.
\end{equation}
\[
= \sum_{n=0}^{\infty} \fr{a^{2n}q^{n^2} (aq;q)_n}{(a^2 q;q)_{2n} (q;q)_n}
\sum_{k=0}^n a^k q^{k^2} {n\brack k}_q.
\]
Now let $a=q$ in (\ref{help-1}) and simplify to obtain
\[
 \fr{1}{(q;q)_{\infty}} \sum_{n=0}^{\infty} (-1)^n (1-q^{2n+1}) q^{\fr{13n^2+11n}{2}}
 =
(1-q) \sum_{n=0}^{\infty} \fr{q^{n^2+2n}}{(q;q)_{2n+2}} \sum_{k=0}^n q^{k^2+k} {n\brack k}_q .
\]
That is
\[
\fr{1}{(q;q)_{\infty}} \sum_{n=0}^{\infty} (-1)^n (1-q^{2n+1}) q^{\fr{13n^2+11n}{2}}
=
\sum_{n=0}^{\infty} \fr{q^{n^2+2n}}{(q^2;q)_{2n+1}} \sum_{k=0}^n q^{k^2+k} {n\brack k}_q,
\]
which completes the proof of the theorem. As to the corollary, simply note that
the left-hand side of previous formula is the generating function of
\[
\sum_{j=0}^n (-1)^j \Big( p \big(n-j(13j+11)/2 \big)-p \big(n-j(13j+15)/2-1 \big) \Big)
\]
and that its right-hand side has nonnegative coefficients.
\section{Proof of Theorem~\ref{thm result-2} and Corollary~\ref{cor result-2}}
\noindent
We need the following formula of Berkovich and Warnaar~\cite[p. 36]{Andrews-Berkovich}
\[
\sum_{k=0}^{\lfloor{n/4\rfloor}}
\fr{(1-a^2 q^{16 k})(a^2, a^2q^4/b;q^8)_k (q^{-n};q)_{4k}}{(1-a^2) (q^8,bq^4;q^8)_k (aq^{n+1};q)_{4k}}
(bq^{4n-2})^k
\]
\[
= \fr{q^n (-q^{-1};q^2)_n (aq;q)_n}{(-q;q)_n (aq;q^2)_n}
\sum_{k\geq 0}\fr{(b^{1/2},-b^{1/2},-aq^4,q^{-2n},q^{2-2n};q^4)_k q^{4k}}
{(q^4, b,aq^2,-q^{3-2n},-q^{5-2n};q^4)_k}
\]
which is equivalent with
\begin{equation}\label{Berk-Warn 1}
\sum_{k=0}^{\lfloor{n/4\rfloor}}
\fr{(1-a^2 q^{16 k})(a^2, a^2q^4/b;q^8)_k (q;q)_n (aq;q)_n q^{{4k\choose 2} -2k} b^k}{(1-a^2) (q^8,bq^4;q^8)_k (q;q)_{n-4k}(aq;q)_{n+4k}}
\end{equation}
\[
= \fr{q^n (-q^{-1};q^2)_n (aq;q)_n}{(-q;q)_n (aq;q^2)_n}
\sum_{k \geq 0}\fr{(b^{1/2},-b^{1/2},-aq^4,q^{-2n},q^{2-2n};q^4)_k q^{4k}}
{(q^4, b,aq^2,-q^{3-2n},-q^{5-2n};q^4)_k}.
\]
We apply Lemma~\ref{lem transf-1} to (\ref{Berk-Warn 1}) as follows.
Let $m=4$, $d= e=1$ and let
\[
A_n= \fr{b^n q^{4n(2n-1)}(1-a^2 q^{16n}) (a^2,a^2q^4/b;q^8)_n}{(1-a^2)(q^8,bq^4;q^8)_n},\quad
D_n = a^{n} q^{n^2},
\]
\[
U_n=\fr{1}{(q;q)_n},\quad V_n=\fr{1}{(a q;q)_n}.
\]
Then on the one hand by~(\ref{Berk-Warn 1}),
\[
\begin{split}
B_n &= \sum_{j=0}^{\lfloor n/4\rfloor} A_j U_{n-4j}V_{n+4j} \\
&= \fr{q^n (-q^{-1};q^2)_n }{(-q;q)_n (q;q)_n (aq;q^2)_n}
\sum_{k\geq 0}\fr{(b^{1/2},-b^{1/2},-aq^4,q^{-2n},q^{2-2n};q^4)_k q^{4k}}
{(q^4, b,aq^2,-q^{3-2n},-q^{5-2n};q^4)_k}.
\end{split}
\]
On the other hand by~(\ref{Bailey 0}),
\[
C_n = \sum_{k=4n}^{\infty} \fr{a^k q^{k^2}}{(q;q)_{k-4n}(aq;q)_{k+4n}}
= \fr{a^{4n} q^{16 n^2}}{(aq;q)_{\infty}}.
\]
Now apply Lemma~\ref{lem transf-1} to these sequences to obtain
from
\[
\sum_{n=0}^{\infty} B_n D_n = \sum_{n=0}^{\infty} A_n C_n
\]
that
\begin{equation}\label{help 2-0}
\sum_{n=0}^{\infty} \fr{a^n q^{n^2+n}(-q^{-1};q^2)_n}{(q^2;q^2)_n (aq;q^2)_n}
\sum_{k\geq 0} \fr{ q^{4k} (b^{1/2},-b^{1/2},-aq^4,q^{-2n}, q^{2-2n};q^4)_k}
{(q^4,b,aq^2,-q^{3-2n},-q^{5-2n};q^4)_k}
\end{equation}
\[
= \fr{1}{(aq;q)_{\infty}} \sum_{n=0}^{\infty} a^{4n} b^n q^{24 n^2 -4n}
\fr{(1-a^2 q^{16n})(a^2,a^2q^4/b;q^8)_n} {(1-a^2) (q^8,bq^4;q^8)_n}
\]
Letting in~(\ref{help 2-0}) $b\to 0$ and then $a=q^4$ we obtain
\[
\fr{1}{(q^5;q)_{\infty}}\sum_{n=0}^{\infty} (-1)^n q^{40 n^2-4n+4} \fr{1- q^{16n+8}}{1-q^8}
\]
\[
=\sum_{n=0}^{\infty}\fr{q^{n^2+5n}(-q^{-1};q^2)_n}{(q^2;q^2)_n (q^5;q^2)_n}
\sum_{k\geq 0} \fr{(-q^8;q^4)_k (q^{-2n};q^2)_{2k}}{(q^4;q^4)_k (q^6;q^4)_k (-q^{3-2n};q^2)_{2k} q^{4k}},
\]
which after some straightforward simplification yields
\[
\fr{1}{(q)_{\infty}}\sum_{n=0}^{\infty} (-1)^n q^{40 n^2-4n+4} (1-q^{16n+8})
\]
\[
= \sum_{n=0}^{\infty}\fr{q^{n^2+5n} (-q^{-1};q^2)_n}{(q;q)_{2n+3}}
\sum_{k=0}^n \fr{q^{ 4k^2+2k-4nk}(-q^4;q^4)_{k+1} (q^2;q^2)_{2n}}{(q^2;q^2)_{2k+1}
(q^2;q^2)_{2n-2k} (-q^{3-2n};q^2)_{2k}},
\]
or equivalently,
\begin{equation}\label{help 2-1}
\fr{1}{(q)_{\infty}}\sum_{n=0}^{\infty} (-1)^n q^{40 n^2-4n+4} (1-q^{16n+8})
\end{equation}
\[
=
\sum_{n=0}^{\infty}\fr{(-q^{-1};q^2)_n}{(q;q)_{2n+3}}
\sum_{k=0}^n \fr{q^{(n-2k)^2 + 5n + 2k}(-q^4;q^4)_{k+1}}{(1-q^{4k+2})(-q^{3-2n};q^2)_{2k}}
{2n \brack 2k}_{q^2}.
\]
This proves the theorem. Concerning the corollary,
multiply both sides of~(\ref{help 2-1}) by $(-q;q)_{\infty}$.
Then by (\ref{generating}) the left-hand side of the obtained formula is the generating function of
\[
\sum_{j=0}^n (-1)^j \Big( \overline{p} \big(n-4j(10j+1)-4 \big)-
\overline{p} \big(n-4j(10j+3)-12 \big) \Big).
\]
Moreover, it is easily seen that
the right-hand side of the obtained formula has nonnegative coefficients.
This completes the proof.
\section{Proof of Theorem~\ref{thm result-3} and Corollary~\ref{cor result-3} }
\noindent
As to this result we shall appeal to the following formula of Berkovich and
Warnaar~\cite[p. 36]{Berkovich-Warnaar}
\[
\sum_{k=0}^{\lfloor{n/4\rfloor}}
\fr{(1-a^4 q^{32 k}) (a^4;q^{16})_k (q^{-n};q)_{4k}}{(1-a^4)(q^{16};q^{16})_k ) (a q^{n+1};q)_{4k}}
(-q^{4n-6})^k
\]
\[
=q^n \fr{(aq^{-1};q^2)_n (aq;q)_n}{(-q;q)_n (aq;q^2)_n} \sum_{k\geq 0}
\fr{(iq^{-2},-iq^{-2},-aq^4, q^{-2n},q^{2-2n};q^4)_k}{(q^4,-q^4,aq^2,q^{3-2n},q^{5-2n};q^4)_k} q^{8k}
\]
which after simplification yields
\begin{equation}\label{Berk-Warn 2}
\sum_{k=0}^{\lfloor{n/4\rfloor}} (-1)^k q^{8k^2-8k} \fr{(1-a^4 q^{32k})(a^4;q^{16})_k}{(1-a^4)
(q^{16};q^{16})_k (q;q)_{n-4k} (aq;q)_{n+4k}}
\end{equation}
\[
= \fr{q^n (-q^{-1};q^2)_n}{(q^2;q^2)_n (aq;q^2)_n} \sum_{k\geq 0}
\fr{q^{8k} (-q^4;q^8)_k (q^{-2n};q^2)_{2k} (-aq^4;q^4)_k}
{(-q^8;q^8)_k (aq^2;q^4)_k (q^{3-2n};q^2)_{2k}}.
\]
Let $m=4$ and $d= e=1$ and let
\[
A_n= \fr{(-1)^n q^{2n(4n-4)}(1-a^4 q^{32n}) (a^4;q^{16})_n}{(1-a^4)(q^{16};q^{16})_n},\quad
D_n = a^{n} q^{n^2},
\]
\[
U_n=\fr{1}{(q;q)_n},\quad V_n=\fr{1}{(a q;q)_n}.
\]
Then by~(\ref{Berk-Warn 2}),
\[
\begin{split}
B_n &= \sum_{j=0}^{\lfloor n/4\rfloor} A_j U_{n-4j}V_{n+4j} \\
&= \fr{q^n (-q^{-1};q^2)_n}{(q^2;q^2)_n (aq;q^2)_n} \sum_{k\geq 0}
\fr{q^{8k} (-q^4;q^8)_k (q^{-2n};q^2)_{2k} (-aq^4;q^4)_k}
{(-q^8;q^8)_k (aq^2;q^4)_k (q^{3-2n};q^2)_{2k}}
\end{split}
\]
and by~(\ref{Bailey 0}),
\[
C_n = \sum_{k=4n}^{\infty} \fr{a^k q^{k^2}}{(q;q)_{k-4n}(aq;q)_{k+4n}}
= \fr{a^{4n} q^{16 n^2}}{(aq;q)_{\infty}}.
\]
Thus by virtue of Lemma~\ref{lem transf-1}
\[
\sum_{n=0}^{\infty} B_n D_n = \sum_{n=0}^{\infty} A_n C_n
\]
and therefore
\[
\sum_{n=0}^{\infty} \fr{a^n q^{n^2+n} (-q^{-1};q^2)_n}{(q^2;q^2)_n (aq;q^2)_n}
\sum_{k\geq 0} \fr{q^{8k} (-q^4;q^8)_k (q^{-2n};q^2)_{2k} (-aq^4;q^4)_k}
{(-q^8;q^8)_k (q^{3-2n};q^2)_{2k} (aq^2;q^4)_k}
\]
\[
= \fr{1}{(aq;q)_{\infty}} \sum_{n=0}^{\infty} (-1)^n a^{4n} q^{24 n^2-8n}\fr{(1-a^4 q^{32n})(a^4;q^{16})_n}
{(1-a^4) (q^{16};q^{16})_n}.
\]
Letting  in the foregoing identity $a=q^4$ and simplifying give
\[
\fr{1}{(q^5;q)_{\infty}} \sum_{n=0}^{\infty} (-1)^n q^{24n^2 +8n} \fr{1-q^{32n+16}}{1-q^8}
\]
\[
= \sum_{n=0}^{\infty} \fr{q^{n^2+5n}(-q^{-1};q^2)_n}{(q^2;q^2)_n (q^5;q^2)_n}
\sum_{k\geq 0} \fr{1+q^{4k+4}}{1+q^4} \fr{q^{4k^2+6k-4nk}(-q^4;q^8)_k (q^2;q^2)_{2n}}
{(q^4;q^4)_k (q^6;q^4)_k (q^2;q^2)_{2n-2k}(q^{3-2n};q^2)_{2k}}
\]
which by rearranging becomes
\[
\fr{1}{(q;q)_{\infty}} \sum_{n=0}^{\infty} (-1)^n q^{8n(3n+1)}(1-q^{32n+16})
\]
\[
= \sum_{n=0}^{\infty}\fr{(-q^{-1};q^2)_n}{(q^2;q^2)_n (q;q^2)_{n+2}}
\sum_{k= 0}^n \fr{q^{(n-2k)^2 + 5n + 6k}(1+q^{4k+4})}{1-q^{4k+2}} {2n \brack 2k}_{q^2}.
\]
which proves the theorem. The corollary follows immediately as the left-hand side of the foregoing formula
is the generating function of
\[
\sum_{j=0}^n (-1)^j \Big( p \big(n-8j(3j+1) \big)-
p \big(n-8j(3j+5)-16 \big) \Big)
\]
and its right-hand side has nonnegative coefficients.
\section{Proof of Theorem~\ref{thm result-4} and Corollary~\ref{cor result-4}}
\noindent
The proof relies on the following relation due to Berkovich and Warnaar~\cite[p. 36]{Berkovich-Warnaar}
\[
\sum_{k=0}^{\lfloor n/6 \rfloor} \fr{ (1-a^2 q^{24 k}) (a^2;q^{12})_k (q^{-n};q)_{6k}}
{(q^{12};q^{12})_k (aq^{n+1};q)_{6k}} (-a^2 q^{6n-3})^k
\]
\[
= \fr{q^n (-q^{-1};q^2)_n (aq;q)_n}{(-q;q)_n (aq;q^2)_n} \sum_{k\geq 0}
\fr{q^{4k} (a^{2/3}, \omega a^{2/3}, \omega^2 a^{2/3}, -aq^4, q^{-2n}, q^{2-2n};q^4)_k}
{(q^4, a, -a, aq^2, -q^{3-2n},-q^{5-2n};q^4)_k}
\]
where $\omega$ is a third root of unity.
Then the previous formula becomes
\begin{equation}\label{Berk-Warn 3}
\sum_{k=0}^{\lfloor{n/6\rfloor}} (-1)^k a^{2k}q^{18k^2-6k} \fr{(1-a^2 q^{24k})(a^2;q^{12})_k}{(1-a^2)
(q^{12};q^{12})_k (q;q)_{n-6k} (aq;q)_{n+6k}}
\end{equation}
\[
= \fr{q^n (-q^{-1};q^2)_n}{(q^2;q^2)_n (aq;q^2)_n} \sum_{k\geq 0}
\fr{q^{4k} (a^{2/3},\omega a^{2/3},\omega^2 a^{2/3};q^4)_k (q^{-2n};q^2)_{2k} (-aq^4;q^4)_k}
{(q^4,a,-a,aq^2;q^4)_k (-q^{3-2n};q^2)_{2k}}.
\]
%Then to apply Lemma~\ref{lem transf-1} we have the following choices which are suggested
%by~(\ref{Berk-Warn 3}).
Let $m=6$, $d= e=1$ and let
\[
A_n= \fr{(-1)^n a^{2n}q^{3n(6n-2)}(1-a^2 q^{24n}) (a^2;q^{12})_n}
{(1-a^2)(q^{12};q^{12})_n},\quad
D_n = a^{n} q^{n^2},
\]
\[
U_n=\fr{1}{(q;q)_n},\quad V_n=\fr{1}{(a q;q)_n}.
\]
Thus by Lemma~\ref{lem transf-1} applied to these items we find
\[
\sum_{n=0}^{\infty}\fr{a^n q^{n^2+n}(-q^{-1};q^2)_n}{(q^2;q^2)_n (aq;q^2)_n} \sum_{k\geq 0}
\fr{q^{4k}(a^{2/3}, \omega a^{2/3}, \omega^2 a^{2/3},-aq^4;q^4)_k (q^{-2n};q^2)_{2k}}
{(q^4,a,-a,aq^2;q^4)_k (-q^{3-2n};q^2)_{2k}}
\]
\[
=\fr{1}{(aq;q)_{\infty}} \sum_{n=0}^{\infty} (-1)^n a^{8n} q^{54 n^2 -6n} \fr{(1-a^2 q^{24n})(a^2;q^{12})_n}
{(1-a^2)(q^{12};q^{12})_n}.
\]
Letting $a=q^6$ in the forgoing formula we find after easy simplification,
including basic facts on $\omega$,
\[
\fr{1}{(q^7;q)_{\infty}} \sum_{n=0}^{\infty} (-1)^n q^{54 n^2 + 42 n} \fr{1-q^{24n+12}}{1-q^{12}}
\]
\[
=\sum_{n=0}^{\infty}\fr{q^{n^2+7n} (-q^{-1};q^2)_n}{(q^2;q^2)_n (q^7;q^2)_n}
\sum_{k\geq 0} \fr{q^{4k}(1+q^4+q^8) (-q^{10};q^4)_k (q^{-2n};q^2)_{2k}}
{(q^{12};q^8)_k (q^8;q^4)_k (-q^{3-2n};q^2)_{2k}},
\]
or equivalently,
\[
\fr{1}{(q^6;q)_{\infty}} \sum_{n=0}^{\infty} (-1)^n q^{54 n^2 + 42 n} (1-q^{24n+12})
\]
\[
=\sum_{n=0}^{\infty}\fr{q^{n^2+7n} (-q^{-1};q^2)_n}{(q^2;q^2)_n (q^7;q^2)_n}
\sum_{k\geq 0} \fr{q^{4k}(1+q^4+q^8) (1+q^{6+4k}) (q^{-2n};q^2)_{2k}}
{(q^{6};q^4)_k (q^8;q^4)_k (-q^{3-2n};q^2)_{2k}},
\]
which yields,
\[
\fr{1}{(q;q)_{\infty}} \sum_{n=0}^{\infty} (-1)^n q^{54 n^2 + 42 n} (1-q^{24n+12})
\]
\[
=\sum_{n=0}^{\infty}\fr{q^{n^2+7n} (-q^{-1};q^2)_n}{(q^2;q^2)_n (q;q^2)_{n+3}}
\sum_{k\geq 0} \fr{q^{4k^2 + 2k -4nk}(1+q^4+q^8) (1+q^{6+4k}) (q^2;q^2)_{2n}}
{(q^{2};q^4)_{k+1} (q^4;q^4)_{k+1} (q^2;q^2)_{2n-2k}(-q^{3-2n};q^2)_{2k}}.
\]
That is,
\[
\fr{1}{(q;q)_{\infty}} \sum_{n=0}^{\infty} (-1)^n q^{54 n^2 + 42 n} (1-q^{24n+12})
\]
\[
=\sum_{n=0}^{\infty}\fr{(-q^{-1};q^2)_n}{(q^2;q^2)_n (q;q^2)_{n+3}}
\sum_{k=0}^n \fr{q^{(n-2k)^2 + 7n + 2k}(1+q^4+q^8) (1+q^{6+4k}) }
{(1-q^{2k+2})(1-q^{2k+4}) (-q^{3-2n};q^2)_{2k}} {2n \brack 2k}_{q^2}
\]
which confirms the theorem. To prove the corollary, first multiply both sides of the previous identity
by $(-q;q)_{\infty}$. Next observe that
the left-hand side of the obtained formula is the generating function of
\[
\sum_{j=0}^n (-1)^j \Big( \overline{p} \big(n-6j(9j+7) \big)-
\overline{p} \big(n-6j(9j+11)-12 \big) \Big)
\]
and the right-hand side of the obtained formula has nonnegative coefficients.
This completes the proof.
\section{Concluding remarks and open questions}
\noindent
Each one of our proofs for Theorems~1-4
relies on substituting specific values for the parameter $a$ which lead to the appropriate formulas. By choosing different values for $a$ one might derive other important results. For instance, upon letting
in the relation~(\ref{help-1}) $a\to 1$ instead of $a=q$, we get after simplification
\begin{equation}\label{concluding-0}
\fr{1}{(q;q)_{\infty}} \sum_{n=0}^{\infty} (-1)^n (1+q^n) q^{\fr{13n^2-n}{2}}
=
\sum_{n=0}^{\infty}\sum_{k=0}^n \fr{q^{n^2+k^2}}{(q;q)_{2n}}  {n\brack k}_q
\end{equation}
or equivalently
\begin{equation}\label{concluding-1}
\fr{1}{(q;q)_{\infty}} \sum_{n=-\infty}^{\infty} (-1)^n q^{\fr{13n^2+n}{2}}
=
\sum_{n=0}^{\infty}\sum_{k=0}^n \fr{q^{n^2+k^2}}{(q;q)_{2n}}  {n\brack k}_q.
\end{equation}
In addition, Andrews~\cite[ formula (5.8) for $k=1$]{Andrews 1984} and
Warnaar~\cite[ Theorem 1.2 for $k=4$]{Warnaar-2001} evaluated the right-hand side
of~(\ref{concluding-1}) as follows
\begin{equation}\label{concluding-2}
\sum_{n=0}^{\infty}\sum_{k=0}^n \fr{q^{n^2+k^2}}{(q;q)_{2n}}  {n\brack k}_q
= \fr{(q^6,q^7,q^{13};q^{13})_{\infty}}{(q;q)_{\infty}}.
\end{equation}
Then by a combination of (\ref{concluding-1}) and (\ref{concluding-2}) we deduce
\begin{equation}\label{concluding-3}
\fr{1}{(q;q)_{\infty}} \sum_{n=-\infty}^{\infty} (-1)^n q^{\fr{13n^2+n}{2}}
= \fr{(q^6,q^7,q^{13};q^{13})_{\infty}}{(q;q)_{\infty}}
\end{equation}
which is also easily obtained by Jacobi's triple product identity~\cite[p. 15]{Gasper-Rahman}
\begin{equation}\label{Jacobi-triple}
(zq,z^{-1}q,q2;q^2)_{\infty} = \sum_{n=-\infty}^{\infty} (-1)^n q^{n^2} z^n
\end{equation}
applied to $z=q^{1/2}$ and $q$ replaced by $q^{13/2}$.
Therefore both~(\ref{concluding-1}) and~(\ref{concluding-3}) imply that the power series
\[
\fr{1}{(q;q)_{\infty}} \sum_{n=-\infty}^{\infty} (-1)^n q^{\fr{13n^2+n}{2}}
\]
has nonnegative coefficients. Besides, it turns out that for any nonnegative odd integer $M$ we find
by~(\ref{Jacobi-triple}) that
\[
\fr{1}{(q;q)_{\infty}} \sum_{n=-\infty}^{\infty} (-1)^n q^{\fr{Mn^2+n}{2}}
= \fr{(q^{(M-1)/2},q^{(M+1)/2},q^{M};q^{M})_{\infty}}{(q;q)_{\infty}}
\]
showing that the left-hand side of this formula has nonnegative integer coefficients. However, we do not how to arrive at this fact using Lemma~\ref{lem transf-1}. This would suggest a general pattern
for our proofs in Theorems~1-4. Obviously, the first challenging part is how to
obtain a terminating identity of the form
\[
\sum_{k=0}^{\lfloor n/m \rfloor} \fr{A_k(a,q)}{(q;q)_{n-mk}(aq;q)_{n+mk}} = B_n (a,q),
\]
where $A_k (a,q)$ and $B_n(a,q)$ are rational functions of $a$ and $q$.
Furthermore, as each one of the alternating sums in Corollaries~1-4 are nonnegative, it is natural to ask
the question whether combinatorial interpretations exist for such sums.
\\ \\
\noindent
{\bf Acknowledgment.} The author thanks Wadim Zudilin for bringing his attention to the
papers by Warnaar~\cite{Warnaar-2001, Warnaar-2003} and the paper by Berokich
and~Warnaar~\cite{Berkovich-Warnaar} along with other papers.
The author is grateful to him and to Ole Warnaar
for reading carefully a previous version of the paper and for fruitful comments and interesting suggestions which
manifestly improved the presentation and quality of this paper.
\end{document}